\documentclass[12pt]{elsarticle}
\usepackage{amsmath, amssymb, amsthm}
\usepackage{graphicx,amssymb,amsfonts,amsbsy,amsmath,amsthm}
\usepackage{mathtools}
\usepackage{hyperref}
\usepackage[utf8]{inputenc}
\usepackage{floatrow}
\usepackage{float}
\usepackage{subfig}
\usepackage{amssymb}
\usepackage{amsthm}
\usepackage{epstopdf}
\usepackage{anysize}
\usepackage{color,soul}
\usepackage{graphicx}
\usepackage{multirow}
\newtheorem{theorem}{Theorem}[section]
\newtheorem{lemma}[theorem]{Lemma}

\newtheorem{remark}[theorem]{Remark}

\DeclareMathOperator{\e}{e}
\DeclareMathOperator{\real}{Re}
\DeclareMathOperator{\imag}{Im}

\journal{...}
\input epsf.sty
\usepackage{lineno}
\begin{document}

\begin{frontmatter}
\title{Stability, bifurcation and control of a predator-prey ecosystem with prey herd behaviour against generalist predator with gestation delay}

%
\author[label1]{Rajesh Ranjan Patra\corref{cor1}}
\ead{razes911@gmail.com}
\author[label1]{Sarit Maitra\corref{cor1}}
\ead{sarit2010.nt@gmail.com}
\author[label2]{Soumen Kundu}

\address[label1]{\corref{cor1}Dept. of Mathematics, NIT Durgapur, Durgapur - 713209, India}
\address[label2]{Department of Mathematics, School of Advanced Sciences, VIT-AP University, Amaravati, Andhra Pradesh - 522237, India}
\cortext[cor1]{Corresponding author}

\begin{abstract}
In this paper, we proposed a population model depicting the dynamics of a prey species showing group defence against a generalist predator. The group defence characteristic is represented by a non-monotonic functional response. We have established the local stability of the model around the co-existent equilibrium solution using a local Lyapunov function. Condition for existence Hopf bifurcation is obtained along with its normal form. Numerical simulations have been done to confirm the obtained analytical results as well as to validate the proposed model. Sensitivity analysis of the parameters is performed using Latin hypercube sampling(LHS)/partial rank correlation coefficient(PRCC). Blow-up in the population is controlled using the Z-type dynamic method.
\end{abstract}


\end{frontmatter}

\section{Introduction}\label{intro}

Conserving biodiversity and managing the resources of a particular ecosystem are points of major concern nowadays. Though there has been human intervention in the field of wildlife management, forestry, fishery etc.\cite{perry1998,christensen1996}, there are processes such as interactions among trophic levels, which also contribute to keeping the natural balance. Predator-prey was mathematically represented by Volterra in the 1920s. Later for modelling population dynamics of different ecosystems, the model has been modified to fit in various environmental and population behaviours such as growth rate, carrying capacity, food availability, mating frequency, fertility rate, predation rate etc. Among them, the predation rate plays a crucial role in the co-existence of both species, which is shown by the functional response in a model.  Numerous articles on observations and applications of different types of functional responses such as Holling type (type I, II, III, IV), Beddington-DeAngelis, Rosenzweig-MacArthur, Crowley-Martin etc. have been studied in past years.

The dynamics of prey may have a significant effect due to predators in trophic systems. Depending upon the type of predation, it has been pointed out by Hanski et al.\cite{hanski1991} that while the specialist predators(mustelids) contribute to the multi-annual cycles of rodent populations of northern Europe, the generalist predators, e.g. foxes, buzzards, cats may be responsible in stabilizing the rodent population in southern localities.

Considering a food chain model of generalist predators, specialist predators and a prey population, Upadhyaya et al.\cite{upadhyaya1998} discussed the existence of chaos in an ecological model. Hassel and Varley introduced a new functional response function to consider the behaviour of grouping in 1969\cite{hassell1969}. After that, there had been numerous research articles showing the behaviour of the functional response in various ecosystems. This depicts the situation when in some interacting populations, predators form group to attack prey to maximize the predation rate. Several research papers have been published on the study of various aspects of such an ecosystem. Hsu et al.\cite{hsu2008} studied the global dynamics of such model. Kim and Baek\cite{kim2013} studied the impulsive effect on such system. In such cases, the functional response is both prey and predator density dependent, where the predators are assumed to move forming a dense colony\cite{cosner1999}. Similarly, grouping behaviour is also witnessed in prey species. There are many biological shreds of evidence of group defence\cite{hoogland1976,patra2021,bi2022} and group vigilance (the many-eyes effect)\cite{siegfried1975} where preys form group against attacking predators to protect their species from extinction by minimizing the rate of predation. 
Many authors use different kinds of functional responses to model the defensive behaviour in prey\cite{ajraldi2011,braza2012,geritz2013,djilali2019}. Among some recent works, Zhang et al.\cite{zhang2021} studied a diffusive predator-prey model with a prey population showing group defence that includes aggression efficiency in the functional response, which could induce instability and bifurcation. Batabyal \& Jana\cite{batabyal2021} used Beddington-DeAngelis functional response for group defence in prey with a modification to the response function by replacing the prey density with the square-root response function. 

The study of delay differential equations are very important as the incorporation of delay in a population model brings out more realistic characteristics of the model\cite{kuang1993}. In a non-delayed predator-prey model, the growth of a species seems instantaneous, but in reality, it is not. In natural habitats, the birth of an individual takes some time which is termed as gestation delay. The characteristic of a model with and without delay may differ completely as the presence of time delay in a system can show much complicated dynamics. It can even alter the stability nature of an equilibrium point\cite{chen2018}. In recent times, there are many studies that involve gestation delay in the predator population. Singh et al.\cite{singh2016} considered the dynamics of disease in prey and gestation delay in generalized predators modelled by the modified Leslie-Gower scheme. Agrawal et al.\cite{agrawal2016} studied a food chain model involving a prey population, specialist predators and generalist predators with gestation delay in the specialist population where the authors studied the global stability analysis of the system. In a predator-prey model with predator species presented by modified Leslie-Gower scheme with gestation delay, Yuan et al.\cite{yuan2015} investigated bifurcation analysis of the system with the Michaelis–Menten type harvesting of prey which is more realistic harvesting from a biological and economic point of view.

Control strategies in a predator-prey system may serve several purposes such as proper utilization of resources, maintaining ecological balance etc.\cite{silveira2005}. Generalist species can survive a wide number of different environmental conditions and may depend upon a wide variety of resources, whereas specialist species can survive a comparatively lower range of environmental conditions and prey upon limited food sources. Monophagous species are dependent upon only a particular type of food source.

The paper is organized as follows. In section \ref{modelsec}, a population model is formulated with certain assumptions. Some preliminary results derived from the model system are discussed in section \ref{prelimsec}. The local stability analysis of the interior equilibrium point in performed in section \ref{lyapsec} using a Lyapunov function. Section \ref{hopfsec} shows the Hopf bifurcation analysis of the model in the presence of delay. In section \ref{dir_st_hopfsec}, the direction and stability of Hopf bifurcation is discussed. Section \ref{numericalsec} contains the numerical simulations. Section \ref{controlsec} includes the case of finite time blow-up in the population and its control. Finally, the paper ends with the conclusion of our work in section \ref{conclusionsec}.
\section{The model} \label{modelsec}
In this paper, we have formulated and investigated a predator-prey model with defensive prey species and generalized predator species, which relies on alternative food sources for its survival where the prey species is the favourite food of the predator species. Hence, the predator dynamics follow the modified Leslie-Gower scheme. The predator growth is due to sexual reproduction. We also assume the presence of gestation delay in the predator species.\\
A logistic prey population is given by
\begin{equation}
\frac{du}{d\delta}=ru\left(1-\frac{u}{K}\right)-v\;g(u)\nonumber\\
\end{equation}
where, $u(\delta)>0$ is the prey population density, $v(\delta)>0$ is the predator population density at any time $\delta$ and $g(u,v)$ is the functional response which shows the defensive behaviour of the prey species. The defensive behaviour is adopted by the preys to reduce the predation rate and hence the ability would be increased when they are in large number. So, as discussed in \cite{patra2021}, the behaviour should be modelled by a non-monotonic functional response with the following characteristics
\begin{eqnarray*}
g:[0,\infty)\mapsto\mathbb{R},\; g \text{ is continuously differentiable},\nonumber\\
g(0)=0,\; g(x)>0 \;\&\; g(x)\leq M \text{ for } \forall x>0,\\
\text{ and } \exists \bar{g}>0\text{ s.t. } g'(x)\begin{cases}
>0, & x<\bar{g},\\
<0, & x>\bar{g}.
\end{cases}
\nonumber
\end{eqnarray*}
as discussed Xiao and Ruan\cite{ruan2001}. The simplified Monod-Haldane functional response, proposed by Sokol and Howell\cite{sokol1981}, fits into the above criteria. The response takes the form
\begin{equation*}
g(u)=\frac{\gamma u}{u^2+\alpha_1}.
\end{equation*}
Let us define $g(u)=\frac{\gamma u}{u^p+\alpha_1}$, where $g(u)$ is same as the simplified Monod-Haldane function for $p=2$. Moreover, $p$ can take other values for which the response satisfies the above requirements. We can notice that, $g(u)$ satisfies the requirements for $p>1$. 

The Leslie-Gower model estimates the carrying capacity of predators species as a quantity directly proportional to the prey density of the environment. The modified Leslie-Gower scheme\cite{aziz2003} incorporates the need for alternating food resources for the predator species when prey density is scarce. Also, the growth of the predator population depends on the number of male and female individuals in the species.\\
Considering all above assumptions, the system can be represented as
\begin{eqnarray}\label{sys0}
\frac{dX}{dT} &=& RX\left(1-\frac{X}{K}\right)-\frac{MXY}{X^p+C},\nonumber\\
\frac{dY}{dT} &=& DY^2-\frac{EY^2}{X+A}.
\end{eqnarray} 
All the parameters are considered to be positive.
\begin{table}[ht]
\small \caption{\sf Biological meaning of  the  parameters for model  (\ref{sys0})} 
\centering
\begin{tabular}{l l l }
\hline
Parameters & Description \\ [0.5ex]
\hline
 $R$ & Intrinsic growth rate of the prey species $X$ \\
 $K$ & Carrying capacity for the preys  \\
 $M$ & Maximum predation rate  \\
 $C$ & The protection provided to the prey population by the environment \\
 $D$ & Reproduction rate of the generalist predator by sexual reproduction  \\
 $E$ & Maximum rate of death of predator population \\
 $A$ & Measures the other food sources available for the predator species \\
  \hline 
 \end{tabular}
\label{Table1}
 \end{table}
  
Introducing gestation delay into the predator equation, assuming a constant gestation period of $\rho$ units, the system can be rewritten as,
\begin{eqnarray}\label{sys01}
\frac{dX}{dT} &=& RX\left(1-\frac{X}{K}\right)-\frac{MXY}{X^p+C},\nonumber\\
\frac{dY}{dT} &=& DY^2-\frac{EY\; Y(T-\rho)}{X(T-\rho)+A}.
\end{eqnarray} 
After using the transformations $x=X/K$, $y=Y/K$ and $t=RT$, the non-dimensional form of the system (\ref{sys01}) becomes,
\begin{eqnarray}\label{sys}
\frac{dx}{dt}&=&x(1-x)-\frac{mxy}{x^p+c},\nonumber\\
\frac{dy}{dt}&=&\left(dy-\frac{ey(t-\tau)}{x(t-\tau)+a}\right)y,
\end{eqnarray}
subject to initial conditions:
\begin{eqnarray}\label{initial_cond}
x(\theta) &=& \phi_1(\theta)>0,\nonumber\\
y(\theta) &=& \phi_2(\theta)>0,\;\theta\in [-\tau,0);\phi_1(0),\phi_2(0)>0,
\end{eqnarray}
where, $\tau=R\rho$, $m=\frac{M}{RK^{p-1}}$, $c=\frac{C}{K^p}$, $d=\frac{DK}{r}$, $e=\frac{E}{R}$ and $a=\frac{A}{K}$.
\section{Positivity} \label{prelimsec}
\noindent Integrating the system (\ref{sys}), we have,
\begin{eqnarray}\label{3}
x(t)&=& x_0 \exp\left[\int_0^t{\left(1-x(\alpha)-\frac{m\,y(\alpha)}{x^p(\alpha)+c}\right)d\alpha}\right],\nonumber\\
y(t)&=& y_0 \exp\left[\int_0^t{\left(dy(\alpha)-\frac{ey(\alpha-\tau)}{x(\alpha-\tau)+a}\right) \; d\alpha}\right],
\end{eqnarray}
where, $x_0=x(0)$ and $y_0=y(0)$.\\

From the above expression, we conclude that, $x(t),\, y(t)>0$ if $x_0,\, y_0>0$. Hence, the dynamics will always stay in $\mathbb{R}_+^2$ if the initial conditions are positive. This shows that the system has positive solutions.

The boundedness of the corresponding non-delayed system is discussed in our earlier work\cite{patra2021}. We present the local stability analysis of the system for the case involving the delayed system.
\section{Local Lyapunov Stability} \label{lyapsec}
In this section, we have used a Lyapunov function to show the local stability of the interior equilibrium point of system (\ref{sys}). The linearised system around the point $(x_*,y_*)$ is given by
\begin{eqnarray}\label{linsys}
\frac{d\bar{x}}{dt} &=& m_{11}\bar{x}+m_{12}\bar{y},\nonumber\\
\frac{d\bar{y}}{dt} &=& m_{21}\bar{x}(t-\tau)+m_{22}\bar{y}-m_{22}\bar{y}(t-\tau),
\end{eqnarray}
where $m_{11}=\frac{mpx_*^p}{(x_*^p+c)^2}-x_*$, \quad $m_{12}=\frac{-mx_*}{x_*^p+c}$, \quad $m_{21}=\frac{ey_*^2}{(x_*+a)^2}$\quad \& \quad $m_{22}=dy_*=\frac{ey_*}{x_*+a}$.\\
Let us define the quantities $P(t)$ and $Q(t)$ such that
\begin{equation}
P(t)=\frac{\bar{x}(t)}{x_*} \quad \& \quad Q(t)=\frac{\bar{y}(t)}{y_*}+\frac{m_{21}}{y_*}\int_{t-\tau}^t\bar{x}(s)ds-\frac{m_{22}}{y_*}\int_{t-\tau}^t\bar{y}(s)ds.
\end{equation} 
Rewriting the system (\ref{linsys}), we have,
\begin{eqnarray}\label{lyap_sys}
\frac{dP}{dt} &=& n_{11}\bar{x}+n_{12}\bar{y},\nonumber\\
\frac{dQ}{dt} &=& n_2\bar{x},
\end{eqnarray}
where $n_{11}=m_{11}/x_*$, $n_{12}=m_{12}/y_*$ \& $n_2=m_{21}/y_*$.\\
Let $\bar{X}(t)=(\bar{x}(t),\bar{y}(t))$ and let
\begin{equation}\label{fun_lxt}
L(\bar{X})(t)=\sum_{i=1}^3{k_iL_i(\bar{X}(t))},
\end{equation}
where $L_i$'s$(i=1,2,3)$ are defined as
\begin{eqnarray}\label{li_vals}
L_1(\bar{X})(t) &=& P^2(t),\nonumber\\
L_2(\bar{X})(t) &=& P(t)\cdot Q(t)+\frac{1}{2}n_2(n_{11}+n_{12})\int_{t-\tau}^t\int_s^t\left[\bar{x}^2(l)+\frac{x_*+a}{y_*}\bar{y}^2(l)\right]\;dl\;ds,\\
L_3(\bar{X})(t) &=& Q^2(t)+n_2^2\int_{t-\tau}^t\int_s^t\left[\bar{x}^2(l)-\frac{x_*+a}{y_*}\bar{y}^2(l)\right]\;dl\;ds,\nonumber
\end{eqnarray}
and the constants $k_i$'s are given by
\begin{eqnarray}\label{ki_vals}
&& k_1=k_2=\frac{2e}{(x_*+a)^2},\\
&& k_3=\frac{2m}{x_*(x_*^p+c)}+\frac{1}{y_*}-\frac{mpx_*^{p-1}}{y_*(x_*^p+c)^2}.
\end{eqnarray}
From equations in (\ref{li_vals}), along the solutions of system (\ref{lyap_sys}), time derivatives of $L_i$'s are calculated to be
\begin{eqnarray}
\frac{dL_1}{dt} &=& \frac{2n_{11}}{x_*}\bar{x}^2+\frac{2n_{12}}{x_*}\bar{x}\bar{y},\\
\frac{dL_2}{dt} &=& \frac{n_{11}}{y_*}\bar{x}\bar{y}+\frac{n_{12}}{y_*}\bar{y}^2+\frac{n_2}{x_*}\bar{x}^2+\frac{1}{2}n_2\tau(n_{11}+n_{12})\left(\bar{x}^2+\frac{1}{y_*}\bar{y}^2\right)+n_{11}n_2\int_{t-\tau}^t\bar{x}(t)\;\bar{x}(s)\;ds\nonumber\\
&& +n_{12}n_2\int_{t-\tau}^t\bar{y}(t)\;\bar{x}(s)\;ds+\frac{x_*+a}{y_*}n_{11}n_2\int_{t-\tau}^t\bar{x}(t)\;\bar{y}(s)\;ds,\nonumber\\
&& +\frac{x_*+a}{y_*}n_{12}n_2\int_{t-\tau}^t\bar{y}(t)\;\bar{y}(s)\;ds,\\
\frac{dL_3}{dt} &=& \frac{2n_2}{y_*}\bar{y}+2n_2^2\int_{t-\tau}^t\bar{x}(t)\;\bar{x}(s)\;ds-\frac{x_*+a}{y_*}2n_2^2\int_{t-\tau}^t\bar{x}(t)\;\bar{y}(s)\;ds+n_2^2\tau\bar{x}^2+\frac{x_*+a}{y_*}n_2^2\tau\bar{y}^2\nonumber\\
&& -n_2^2\int_{t-\tau}^t\bar{x}^2(s)\;ds-\frac{x_*+a}{y_*}n_2^2\int_{t-\tau}^t\bar{y}^2(s)\;ds.
\end{eqnarray}
Using the inequality $ab\leq\frac{1}{2}(a^2+b^2)$ \& (\ref{fun_lxt}), we have
\begin{eqnarray}\label{dldt}
\frac{dL}{dt}(\bar{X})(t)&\leq & U_1\bar{x}^2+U_2\bar{y}^2,
\end{eqnarray}
where
\begin{eqnarray*}
U_1 &=& k_1\frac{2n_{11}}{x_*}+k_2\left(\frac{n_{11}n_{2}}{2}+\frac{n_{11}}{2}\frac{e}{x_*+a}+\frac{n_2}{x_*}+\frac{n_{11}n_2\tau}{2}\right)+k_3\left(\frac{n_{12}n_2\tau}{2}+n_2^2+n_2^2\tau+\frac{e}{x_*+a}n_2\right),\nonumber\\
U_2 &=& k_2\left(\frac{n_{12}}{y_*}+\frac{a_{12}a_2}{2}+\frac{a_{12}}{2}\frac{e}{x_*+a}+\frac{a_{11}\tau}{2}\frac{e}{x_*+a}+\frac{a_{12}\tau}{2}\frac{e}{x_*+a}\right)+k_3\frac{e}{x_*+a}n_2\tau.
\end{eqnarray*}
\begin{theorem}\label{lyap_thm}
The equilibrium point $(x_*,y_*)$ is locally asymptotically stable if the quantities $U_1$ and $U_2$ are both negative. 
\end{theorem}
\begin{proof}
Going by the steps as performed in \cite{kundu2016}, we can show that, the function, $L(\bar{X})(t)$, is a Lyapunov function if the delay value $\tau$ satisfies $U_1,U_2<0$, provided $k_i$'s are positive.
Hence, the interior equilibrium point is locally asymptotically stable.
\end{proof}
The following figures depict the plot of $U_1$ and $U_2$ for different values of $\tau$. The first figure shows that both the quantities $U_1$ and $U_2$ are negative where the internal equilibrium point is locally asymptotically stable. The second figure shows that at least one of the quantities $U_1$ and $U_2$ are positive for a range $\tau$ values.
\begin{figure}[H]
 \subfloat[Both $U_1$ and $U_2$ are -ve\label{subfig-7}]{%
  \includegraphics[width=0.45\textwidth, height=6cm]{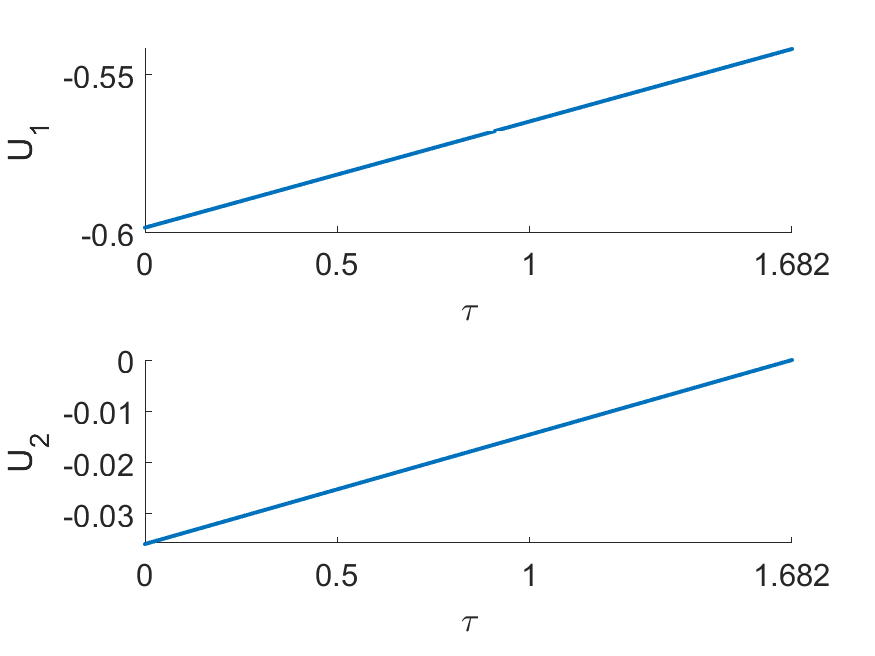}}
 \subfloat[At least one of $U_1$ and $U_2$ is +ve\label{subfig-7}]{%
  \includegraphics[width=0.45\textwidth, height=6cm]{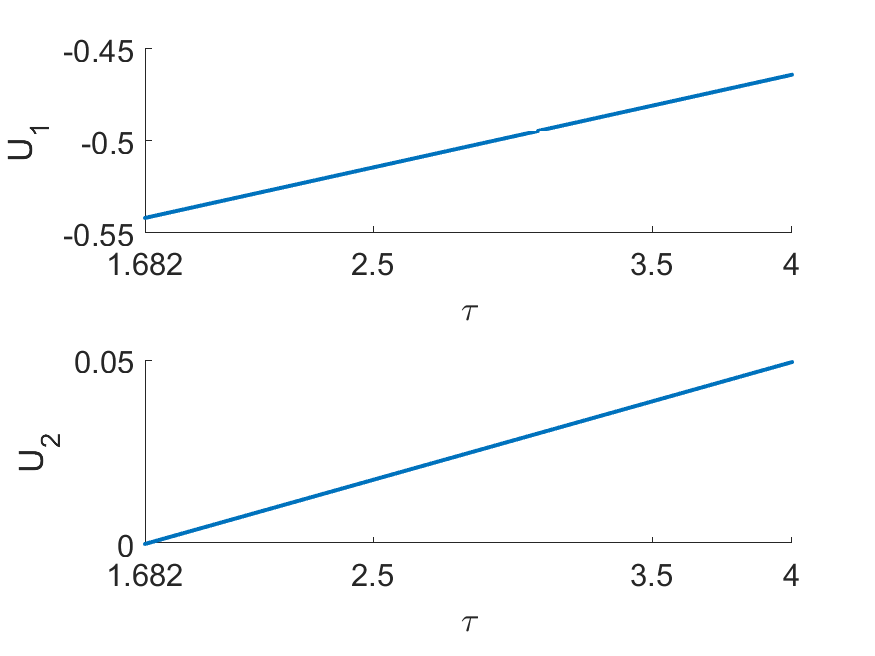}}
 \caption{Plot of $U_1$ and $U_2$ for a range of values of the delay parameter$(\tau)$.}
 \label{bb}
\end{figure}
\section{Hopf Bifurcation}\label{hopfsec}
In a system of differential equations, Hopf bifurcation happens when the complex conjugate set of eigenvalues of a linearised system become purely imaginary at a fixed point. Hopf-bifurcation occurs at a point where a system changes state from stable to unstable, i.e. it is a local bifurcation in which a fixed point of a dynamical system loses stability. In this section, we have studied Hopf bifurcation for system (\ref{sys}) for equilibrium point $(x_*,y_*)$. From the linearised system (\ref{linsys}), the characteristic equation of the system for $(x_*,y_*)$ is given by,
\begin{equation}\label{ch_eqn}
\lambda^2-(m_{11}+m_{22})\lambda+m_{11}m_{22}+(m_{22}\lambda-m_{11}m_{22}-m_{12}m_{21})\e^{-\lambda\tau}=0.
\end{equation}
For $\tau=0$, equation (\ref{ch_eqn}) changes to
\begin{equation}\label{ch_eqn_nd}
\lambda^2-m_{11}\lambda-m_{12}m_{21}=0,
\end{equation}
which has all roots with negative real parts if $m_{11}<0$ as $m_{12}m_{21}<0$.\\
For $\tau>0$, let us substitute $\lambda=i\omega$, where $\omega>0$ \& $i=\sqrt{-1}$, then from the real and imaginary part of (\ref{ch_eqn}), we have,
\begin{eqnarray}\label{re_im_eqns}
-(m_{11}m_{22}+m_{12}m_{21})\cos{\omega\tau}+m_{22}\omega\;\sin{\omega\tau} &=& \omega^2-m_{11}m_{22},\nonumber\\
-m_{22}\omega\;\cos{\omega\tau}-(m_{11}m_{22}+m_{12}m_{21})\;\sin{\omega\tau} &=& \omega(m_{11}+m_{22}).
\end{eqnarray}
Squaring both sides of equations in (\ref{re_im_eqns}) and adding them,
\begin{equation}\label{omega_eqn}
\omega^4+m_{11}^2\omega^2+m_{11}^2m_{22}^2-(m_{11}m_{22}+m_{12}m_{21})^2=0.
\end{equation}
Equation (\ref{omega_eqn}) has a positive root $\omega_0$ if
\begin{eqnarray}\label{omega_cond}
& m_{11}^2m_{22}^2-(m_{11}m_{22}+m_{12}m_{21})^2=-m_{12}m_{21}\left(2m_{11}m_{22}+m_{12}m_{21}\right)<0,
\end{eqnarray}
which holds for $m_{11}<0$ as we have $m_{12}<0$, $m_{21}>0$ \& $m_{22}>0$.\\
Eliminating $\sin{\omega\tau}$ from (\ref{re_im_eqns}),
\begin{eqnarray}\label{cos_omega_tau}
\cos{\omega\tau}=\frac{m_{22}(m_{11}+m_{22})\omega^2-(\omega^2-m_{11}m_{22})(m_{11}m_{22}+m_{12}m_{21})}{m_{22}^2\omega^2+(m_{11}m_{22}+m_{12}m_{21})^2}.
\end{eqnarray}
With $\omega=\omega_0$ in (\ref{cos_omega_tau}), for $k=0,1,2,...$,
\begin{eqnarray}\label{tau_n}
\tau_k=\frac{1}{\omega_0}\arccos{\left[\frac{m_{22}(m_{11}+m_{22})\omega^2-(\omega^2-m_{11}m_{22})(m_{11}m_{22}+m_{12}m_{21})}{m_{22}^2\omega^2+(m_{11}m_{22}+m_{12}m_{21})^2}\right]}+\frac{2\pi k}{\omega_0}.
\end{eqnarray}
By Lemma \ref{butlers_lemma}, stability switching occurs at $\tau=\tau_0$. Differentiating (\ref{ch_eqn}) w.r.t. $\tau$, where $\lambda=\lambda(\tau)$, we get,
\begin{eqnarray}\label{d_lambda_d_tau_inv}
\left(\frac{d\lambda}{d\tau}\right)^{-1} &=& \frac{m_{22}\e^{-\lambda\tau}+2\lambda-m_{11}-m_{22}}{(m_{22}\lambda+m_{11}m_{22}+m_{12}m_{21})\lambda\e^{-\lambda\tau}}-\frac{\tau}{\lambda},
\end{eqnarray}
Now, putting $\tau=\tau_0$ and $\lambda=i\omega_0$,
\begin{equation}\label{re_d_lambda_d_tau}
\real\left[\left(\frac{d\lambda}{d\tau}\right)^{-1}_{\tau=\tau_0, \lambda=i\omega_0}\right]=\frac{2\omega_0^2+m_{11}^2}{m_{22}\omega_0^2+(m_{11}m_{22}+m_{12}m_{21})^2},
\end{equation}
which leads us to
\begin{equation}\label{trans_sign}
\text{Sign}\left[\real\left(\frac{d\lambda}{d\tau}\right)^{-1}_{\tau=\tau_0,\lambda=i\omega_0}\right]>0.
\end{equation}
Hence,
\begin{equation}\label{trans_eqn}
\frac{d}{d\tau}\real(\lambda)>0.
\end{equation}
Therefore, the transversality condition holds which implies that the system undergoes a Hopf bifurcation at $(x_*,y_*)$ as $\tau$ crosses $\tau_0$. 

\begin{lemma}[G. J. Butler's Lemma\cite{freedman1983}]\label{butlers_lemma}
Let $m_{11}<0$. Then the real parts of the solutions of (\ref{ch_eqn}) are negative for all $\tau<\tau_0$, where $\tau_0>0$ is the smallest value for which there is a solution to equation (\ref{ch_eqn}) with real part zero.
\end{lemma}

Using Lemma \ref{butlers_lemma}, if $m_{11}<0$, the equilibrium point is locally stable for delay value $\tau<\tau_0$.
\section{Direction and stability of Hopf bifurcation}\label{dir_st_hopfsec}
In this section, we will use the normal form and the center manifold theory\cite{hassard1981} to analyze the stability, direction and period of the periodic solutions around the equilibrium point $(x^*,y^*)$ originating as $\tau$ crosses the critical value $\tau_0$. In section \ref{hopfsec}, we have obtained the condition for occurrence of Hopf bifurcation in system (\ref{sys}). So, as derived earlier, Hopf bifurcation occurs at $(x^*,y^*)$ for $\tau=\tau_0$ and the eigenvalues of the Jacobian of the system at that point are $\pm i\omega_0$, and we have $$\left(\frac{d\lambda}{d\tau}\right)_{\tau=\tau_0}\neq 0.$$

For convenience, let $\tau=\tau_0+\mu$, then $\mu=0$ is the Hopf bifurcation value of (\ref{sys}). Let $u_1(t)=x(t)-x^*$ \& $u_2(t)=y(t)-y^*$, then letting $x(t)=u_1(\tau t)$ \& $y(t)=u_2(\tau t)$, the system converts into a functional differential equation in $C[-1,0],R^2)$ as
\begin{equation}\label{fde}
\frac{du}{dt}=L_\mu (u_t)+f(\mu,u_t),
\end{equation}
where $u(t)=(x(t),y(t)^T\in R^2$, $u_t(\theta)=u(t+\theta)$, $\theta\in [-1,0]$, and $L_\mu : C\rightarrow R$ \& $f:R \times C\rightarrow R$ are given by
\begin{equation}\label{l_mu_1}
L_\mu\phi = (\tau_0+\mu)B\phi(0)+(\tau_0+\mu)G\phi(-1),
\end{equation}
where
\begin{equation*}
B= \begin{pmatrix}
m_{11} & m_{12}\\
0 & m_{22}
\end{pmatrix},
\quad
G=\begin{pmatrix}
0 & 0\\
m_{21} & -m_{22}
\end{pmatrix}
\quad \& \quad
f(\mu,\phi)=\begin{pmatrix}
f_1(\mu,\phi)\\
f_2(\mu,\phi)
\end{pmatrix},
\end{equation*}
where
\begin{eqnarray}\label{f1_f2}
f_1(\mu,\phi) &=& (\tau_0+\mu)[m_{13}\phi_1^2(0)+m_{14}\phi_1(0)\phi_2(0)+m_{15}\phi^3_1(0)+m_{16}\phi_1^2(0)\phi_2(0)+\text{h.o.t.}],\nonumber\\
f_2(\mu,\phi) &=& (\tau_0+\mu)[m_{23}\phi_2^2(0)+m_{24}\phi_1^2(-1)+m_{25}\phi_1(-1)\phi_2(0)+m_{26}\phi_1(-1)\phi_2(-1)\\
& & +m_{27}\phi_2(0)\phi_2(-1)+m_{28}\phi_1^2(-1)\phi_2(0)+m_{29}\phi_1^3(-1)+m_{30}\phi_1^2(-1)\phi_2(-1)+\text{h.o.t.}],\nonumber
\end{eqnarray}
are obtained using the similar procedure as done in \cite{meng2011} \& \cite{xu2014}, and the coefficients $m_{jk}$'s are mentioned below.
\begin{eqnarray*}
m_{13} &=& \frac{mpx^{*^{p-1}}y^*}{2(x^{*^p}+c)^2}-1-\frac{mp^2x^{*^{p-1}}y^*(x^{*^p}-c)}{2(x^{*^p}+c)^3},\\
m_{14} &=& \frac{mpx^{*^p}}{2(x^{*^p}+c)^2}-\frac{m}{2(x^{*^p}+c)},\\
m_{15} &=& \frac{mp^2x^{*^{p-2}}y^*((p+1)x^{*^{2p}}-4cpx^{*^{2p}}+(p-1)c^2)}{6(x^{*^p}+c)^4}-\frac{mp^2x^{*^{p-2}}y^*((p+1)x^{*^p}-pc+c)}{6(x^{*^p}+c)^3},\\
m_{16} &=& \frac{mpx^{*^{p-1}}}{6(x^{*^p}+c)^2}-\frac{mp^2x^{*^{p-1}}(x^{*^p}-c)}{6(x^{*^p}+c)^3},\\
m_{23} &=& d,\\
m_{24} &=& -\frac{ey^{*^2}}{(x^*+a)^3},\\
m_{25} &=& \frac{ey^*}{2(x^*+a)^2}=m_{26},\\
m_{27} &=& -\frac{e}{x^*+a},\\
m_{28} &=& -\frac{ey^*}{3(x^*+a)^3}=m_{30},\\
m_{29} &=& \frac{ey^{*^2}}{(x^*+a)^4}.
\end{eqnarray*}
By Reisz representation theorem, there exists a function $\eta(\theta,\mu)$ whose components are of bounded variation for $\theta\in [-1,0]$ such that
\begin{equation}\label{l_mu_2}
L_{\mu}\phi=\int_{-1}^0\phi(\theta)\;d\eta(\theta,\mu),\quad \text{for }\phi\in C.
\end{equation}
In view of equation (\ref{l_mu_1}), we can choose
\begin{equation}\label{eta}
\eta(\theta,\mu)=(\tau_0+\mu)[B\;\delta(\theta)+G\;\delta(\theta+1),
\end{equation}
where $\delta$ is the Dirac delta function.\\
For $\phi\in C([-1,0],R^2)$, define
\begin{equation}\label{a_mu}
A(\mu)\phi=\begin{cases}
\frac{d\phi(\theta)}{d\theta}, & \theta\in [-1,0),\\
\int_{-1}^0\phi(s)\;d\eta(s,\mu), & \theta=0,
\end{cases}
\end{equation}
and
\begin{equation}\label{r_mu}
R(\mu)\phi=\begin{cases}
0, & \theta\in [-1,0),\\
f(\mu,\phi), & \theta=0.
\end{cases}
\end{equation}
Then (\ref{fde}) becomes
\begin{equation}\label{u_dot}
\dot{u}_t=A(\mu)u_t+R(\mu)u_t.
\end{equation}
For $\psi\in C([0,1],(R^2)^*)$, define
\begin{equation}\label{a_star_psi}
A^*\psi(s)=\begin{cases}
-\frac{d\psi(s)}{ds}, & s\in (0,1],\\
\int_{-1}^0\psi(-t)\; d\eta^T(t,0), & s=0.
\end{cases}
\end{equation}
Define the bilinear form, for $\phi\in C([-1,0],R^2)$ \& $\psi\in C([0,1],(R^2)^*)$, as
\begin{equation}\label{inner_prod}
\langle \psi,\phi \rangle =\bar{\psi}(0)\;\phi(0)-\int_{-1}^0 \int_{\xi=0}^\theta \bar{\psi}(\xi-\theta)\;\phi(\xi)\;d\eta(\theta)\;d\xi,
\end{equation}
where $ \eta(\theta)=\eta(\theta,0)$.\\
Then $A$ \& $A^*$ are adjoint operators. So, $\pm i\omega_0\tau_0$ are the eigenvalues of the matrices $A$ \& $A^*$. We need to compute the eigenvectors of $A$ \& $A^*$corresponding to eigenvalues $i\omega_0\tau-0$ \& $-i\omega_0\tau_0$, respectively.\\
Suppose $q(\theta)=(\alpha_1,\alpha_2)^T\e^{i\omega_0\tau_0\theta}$ be the eigenvector corresponding to the eigenvalues $i\omega_0\tau_0$, then
\begin{equation}\label{a_q_theta}
Aq(\theta)=i\omega_0\tau_o\;q(\theta),
\end{equation}
which, using (\ref{l_mu_1}), (\ref{l_mu_2}) \& (\ref{a_mu}), for $\theta=0$ gives
\begin{equation*}
(\tau_0+\mu)[B\;q(0)+G\; q(-1)]=i\omega_0\tau_0\;q(0),
\end{equation*}
which further gives
\begin{equation*}
\tau_0\begin{pmatrix}
i\omega_0-m_{11} & -m_{12}\\
-m_{21}\e^{-i\omega_0\tau_0} & i\omega_0-m_{22}+m_{22}\e^{-i\omega_0\tau_0}
\end{pmatrix} \begin{pmatrix}
\alpha_1\\
\alpha_2
\end{pmatrix} =\begin{pmatrix}
0\\
0
\end{pmatrix}.
\end{equation*}
Choosing $\alpha_1=1$, from the above expression, we get,
\begin{equation}\label{gamma}
\alpha_2=\frac{i\omega_0-m_{11}}{m_{12}}=\gamma\text{ (let)}.
\end{equation}
So, $q(\theta)=(1,\gamma)^T\e^{i\omega_0\tau_0\theta},\; \theta\in [-1,0].$\\
Similarly, let $q^*(s)=D(1,\gamma^*)^T\e^{i\omega_0\tau_0 s}$ be the eigenvector of $A^*$ corresponding to the eigenvalue $-i\omega_0\tau_0$, then
$$ A^*q^*(s)=-i\omega_0\tau_0q^*(s),$$
$$\text{i.e., }\tau_0\begin{pmatrix}
i\omega_0+m_{11} & m_{21}\e^{-i\omega_0\tau_0}\\
m_{12} & i\omega_0+m_{22}-m_{22}\e^{-i\omega_0\tau_0}
\end{pmatrix} \begin{pmatrix}
1\\
\gamma^*
\end{pmatrix}=\begin{pmatrix}
0\\
0
\end{pmatrix} $$
which gives
\begin{equation}\label{gamma_star}
\gamma^*=-\frac{i\omega_0+m_{11}}{m_{21}\e^{-i\omega_0\tau_0}}.
\end{equation}
So, $q^*(s)=D(1,\gamma^*)^T\e^{i\omega_0\tau_0s},\; s\in [0,1].$\\
Using normalization condition, we have the relation, $\langle q^*(s),q(\theta)\rangle=1$, by which we can find the value of $D$. Using (\ref{inner_prod}),
\begin{eqnarray*}
\langle q^*(s),q(\theta)\rangle &=& \bar{D}(1+\gamma\gamma^*)-\bar{D}\int_{\theta=-1}^0(1,\bar{\gamma^*})\theta\e^{i\omega_0\tau_0\theta}d\eta(\theta)\times (1,\gamma)^T\\
&=& \bar{D}(1+\gamma\bar{\gamma^*})+\bar{D}\tau_0\e^{-i\omega_0\tau_0}(m_{21}-m_{22}\gamma)\bar{\gamma^*},
\end{eqnarray*}
then
\begin{equation}\label{d_bar}
\bar{D}=\frac{1}{1+\gamma\bar{\gamma^*}+\tau_0\e^{-i\omega_0\tau_0}(m_{21}-m_{22}\gamma)\bar{\gamma^*}}.
\end{equation}
Now, we compute the coordinates to describe the center manifold $C_0$ at $\mu=0$. Let $u_t$ be the solution of (\ref{sys}) when $\mu=0$. Define  
\begin{equation}\label{zt_wt}
z(t)=\langle q^*,u_t\rangle ,\quad W(t,\theta)=u_t(\theta)-2\;\real\{z(t)\;q(\theta)\}.
\end{equation}
On the center manifold $C_0$, we have
\begin{equation}\label{w_series}
W(t,\theta)=W(z(t),\bar{z}(t),\theta)=W_{20}(\theta)\frac{z^2}{2}+W_{11}(\theta)z\bar{z}+W_{02}(\theta)\frac{\bar{z}^2}{2}+W_{30}(\theta)\frac{z^3}{6}+\cdots ,
\end{equation}
and $z$ \& $\bar{z}$ are coordinates for center manifold $C_0$ in the direction of $q^*$ \& $\bar{q^*}$. Note that, $W$ is also real if $u_t$ is real, we consider only real solutions. For solutions $u_t\in C_0$ of (\ref{sys}),
\begin{eqnarray}\label{z_dot_1}
\dot{z}(t) &=& \langle q^*(s),A(0)u_t\rangle+\langle q^*(s),R(0)u_t\rangle\nonumber\\
\implies \dot{z}(t) &=& i\omega_0\tau_0z(t)+\bar{q^*}(0)f(0,u_t(\theta)) 
\end{eqnarray}
which becomes
\begin{equation}\label{z_dot_2}
\dot{z}(t)=i\omega_0\tau_0z(t)+g(z,\bar{z})
\end{equation}
by letting $g(z,\bar{z})=\bar{q^*}(0)f(0,u_t(\theta))=\bar{q^*}(0)f(z(t),\bar{z}(t))$, where
\begin{equation}\label{g_fun}
g(z,\bar{z})=g_{20}(\theta)\frac{z^2}{2}+g_{11}z\bar{z}+g_{02}(\theta)\frac{\bar{z}^2}{2}+g_{21}(\theta)\frac{z^2\bar{z}}{2}+\cdots .
\end{equation}
Hence,
\begin{equation}\label{g_fun_1}
g(z,\bar{z})=\bar{q^*}(0)f(0,u_t(\theta))=\bar{D}(1,\gamma^*)(f_1(0,u_t),f_2(0,u_t))^T,
\end{equation}
where $f_1(0,u_t)$ \& $f_2(0,u_t)$ are obtained from (\ref{f1_f2}).\\
Noticing that $u_t(\theta)=(x_t(\theta),y_t(\theta))^T=W(t,\theta)+zq(\theta)+\bar{z}\bar{q}(\theta)$ \& $q(\theta)=(1,\gamma)^T\e^{i\omega_0\tau_0}$, we have
\begin{eqnarray*}
& x_t(0)= z+\bar{z}+W^{(1)}(0),\\
& y_t(0) = z+\bar{z}+W^{(2)}(0),\\
& x_t(-1) = \e^{-i\omega_0\tau_0}z+\e^{i\omega_0\tau_0}\bar{z}+W^{(1)}(-1),\\
& y_t(-1) = \gamma\e^{-i\omega_0\tau_0}z+\bar{\gamma^*}\e^{i\omega_0\tau_0}\bar{z}+W^{(2)}(-1),
\end{eqnarray*}
where $$W^{(k)}(\theta)=W_{20}^{(k)}(\theta)\frac{z^2}{2}+W_{11}^{(k)}(\theta)z\bar{z}+W_{02}^{(k)}(\theta)\frac{\bar{z}^2}{2}+W_{21}^{(k)}(\theta)\frac{z^2\bar{z}}{2}+\cdots,\quad k=1,2 .$$
From (\ref{g_fun_1}),
\begin{eqnarray*}
g(z,\bar{z}) &=& \bar{D}\tau_0[m_{13}x_t^2(0)+m_{14}x_t(0)y_t(0)+m_{15}x_t^3(0)+m_{16}x_t^2(0)y_t(0)+\cdots]\\
& & +\bar{D}\gamma^*\tau_0[m_{23}y_t^2(0)+m_{24}x_t^2(-1)+m_{25}x_t(-1)y_t(0)+m_{26}x_t(-1)y_t(-1)\\
& & +m_{27}y_t(0)y_t(-1)+m_{28}x_t^2(-1)y_t(0)+m_{29}x_t^3(-1)+m_{30}x_t^2(-1)y_t(-1)+\cdots]
\end{eqnarray*}
which becomes
\begin{eqnarray}\label{g_fun_main}
g(z,\bar{z}) &=& \bar{D}\tau_0\bigl[
\bigl\{m_{13}+m_{14}\gamma+m_{23}\gamma^2\gamma^*+m_{24}\gamma^*\e^{-2i\omega_0\tau_0}+m_{25}\gamma\gamma^*\e^{-i\omega_0\tau_0}+m_{26}\gamma\gamma^*\e^{-2i\omega_0\tau_0}\nonumber\\
& & +m_{27}\gamma^2\gamma^*\e^{-i\omega_0\tau_0}\bigr\}\times z^2+\bigl\{
m_{13}+m_{14}\bar{\gamma}+m_{23}\bar{\gamma}^2\gamma^*+m_{24}\gamma^*\e^{2i\omega_0\tau_0}+m_{25}\bar{\gamma}\gamma^*\e^{i\omega_0\tau_0}\nonumber\\
& & +m_{26}\bar{\gamma}\gamma^*\e^{2i\omega_0\tau_0}+m_{27}\bar{\gamma}^2\gamma^*\e^{i\omega_0\tau_0}\bigr\}\times\bar{z}^2
+\bigl\{2m_{13}+2m_{14}\real (\gamma)+2m_{23}\gamma\bar{\gamma}\gamma^*+2m_{24}\gamma^*\nonumber\\
& & +m_{25}\gamma^*(\gamma\e^{i\omega_0\tau_0}+\bar{\gamma}\e^{-i\omega_0\tau_0})2m_{26}\gamma^*\real (\gamma)+m_{27}\gamma\bar{\gamma}\gamma^*(\e^{i\omega_0\tau_0}+\e^{-i\omega_0\tau_0})\bigr\}\times z\bar{z}\nonumber\\
& & +\bigr\{
2m_{13}\bigl(W_{20}^{(1)}(0)+2W_{11}^{(1)}(0)\bigr)+m_{14}\bigl(W_{20}^{(1)}(0)\bar{\gamma}+W_{20}^{(2)}(0)+2W_{11}^{(1)}(0)\gamma+2W_{11}^{(2)}(0)\bigr)+6m_{15}\nonumber\\
& & +2m_{16}(\bar{\gamma}+2\gamma)+2m_{23}\gamma^*\bigl(W_{20}^{(2)}(0)+2W_{11}^{(2)}(0)\bigr)+2m_{24}\gamma^*\bigl(W_{20}^{(1)}(-1)+2W_{11}^{(1)}(-1)\e^{-i\omega_0\tau_0}\bigr)\nonumber\\
& & +m_{25}\gamma^*\bigl(W_{20}^{(1)}(-1)\bar{\gamma}+W_{20}^{(2)}(0)\e^{i\omega_0\tau_0}+2W_{11}^{(1)}(-1)\gamma+2W_{11}^{(2)}(0)\e^{-i\omega_0\tau_0}\bigr)\nonumber\\
& & +m_{26}\gamma^*\bigl(W_{20}^{(1)}(-1)\bar{\gamma}\e^{i\omega_0\tau_0}+W_{20}^{(2)}(-1)\e^{i\omega_0\tau_0}+2W_{11}^{(1)}(-1)\gamma\e^{-i\omega_0\tau_0}+2W_{11}^{(2)}(-1)\e^{-i\omega_0\tau_0}\bigr)\nonumber\\
& & +m_{27}\gamma^*\bigl(W_{20}^{(2)}(0)\bar{\gamma}\e^{i\omega_0\tau_0}+W_{20}^{(2)}(-1)\gamma+W_{11}^{(2)}(0)\gamma\e^{i\omega_0\tau_0}+W_{11}^{(2)}(-1)\gamma\bigr)\nonumber\\
& & +2m_{28}\gamma^*\bigl(\bar{\gamma}\e^{-2i\omega_0\tau_0}+2\gamma\bigr)+6m_{29}\gamma^*\e^{-i\omega_0\tau_0}+2m_{30}\gamma^*(\bar{\gamma}+2\gamma)\e^{-i\omega_0\tau_0}\}\times\frac{z\bar{z}}{2} \bigr]
\end{eqnarray}
To compute $g_{21}$, we need to compute $W_{20}^{(k)}(\theta)$ \& $W_{11}^{(k)}(\theta)$, for $k=1,2$ \& $\theta=0,-1$.
From (\ref{u_dot}), (\ref{zt_wt}) \& (\ref{z_dot_1}),
\begin{equation*}
\dot{W} = \begin{cases}
AW-2\real\{\bar{q^*}(0)f_0q(\theta)\}, & \theta\in [-1,0),\\
AW-2\real\{\bar{q^*}(0)f_0q(\theta)\}, & \theta=0,
\end{cases}
\end{equation*}
\begin{equation}\label{w_dot}
\text{i.e.,}\quad \dot{W}= AW+H(z,\bar{z},\theta),
\end{equation}
where
\begin{eqnarray}\label{h_series}
H(z,\bar{z},\theta)=H_{20}\frac{z^2}{2}+H_{11}(\theta)z\bar{z}+H_{02}(\theta)\frac{\bar{z}^2}{2}+\cdots .
\end{eqnarray}
From (\ref{w_series}), we have $\dot{W}=W_z\dot{z}+W_{\bar{z}}\dot{\bar{z}}$. Then comparing coefficients with that of (\ref{w_dot}),
\begin{eqnarray}
(A-2i\omega_0\tau_0)W_{20} &=& -H_{20}(\theta),\label{w_h_reln_1}\\
AW_{11} &=& -H_{11}(\theta)\label{w_h_reln_2}.
\end{eqnarray}
For $\theta\in[-1,0)$, we have
\begin{eqnarray}\label{h_series_g}
H(z,\bar{z},\theta) &=& -\bigl(g_{20}q(\theta)+\bar{g}_{02}\bar{q}(\theta)\bigr)\frac{z^2}{2}-\bigl(g_{11}q(\theta)+\bar{g}_{11}\bar{q}(\theta)\bigr)z\bar{z}+\cdots .
\end{eqnarray}
Comparing its coefficients with (\ref{h_series}),
\begin{eqnarray}
H_{20}(\theta) &=& -g_{20}q(\theta)-\bar{g}_{02}\bar{q}(\theta),\label{h_g_reln_1}\\
H_{11}(\theta) &=& -g_{11}q(\theta)-\bar{g}_{11}\bar{q}(\theta).\label{h_g_reln_2}
\end{eqnarray}
From (\ref{w_h_reln_1}), (\ref{h_g_reln_1}) \& from the definition of $A$,
\begin{equation*}\label{w20_dot}
\dot{W}_{20}(\theta)=2i\omega_0\tau_0W_{20}+g_{20}q(0)\e^{i\omega_0\tau_0\theta}+\bar{g}_{02}\e^{-i\omega_0\tau_0\theta},
\end{equation*}
which gives
\begin{equation}\label{w20_soln}
W_{20}(\theta)=\frac{ig_{20}}{\omega_0\tau_0}q(0)\e^{i\omega_0\tau_0}+\frac{i\bar{g}_{02}}{3\omega_0\tau_0}\bar{q}(0)\e^{-i\omega_0\tau_0}+E_1\e^{2i\omega_0\tau_0}.
\end{equation}
Similarly, from (\ref{w_h_reln_2}), (\ref{h_g_reln_2}) \& the definition of $A$,
$$\dot{W}_{11}(\theta)=g_{11}q(\theta)+\bar{g}_{11}\bar{q}(\theta),$$ that gives
\begin{equation}\label{w11_soln}
W_{11}(\theta)=-\frac{ig_{11}}{\omega_0\tau_0}q(0)\e^{i\omega_0\tau_0}+\frac{i\bar{g}_{11}}{\omega_0\tau_0}\bar{q}(0)\e^{-i\omega_0\tau_0}+E_2,
\end{equation}
where, $E_1=\bigl(E_1^{(1)},E_1^{(2)}\bigr)$ \& $E_2=\bigl(E_2^{(1)},E_2^{(2)}\bigr)$ are constant vectors in $R^2$ \& to be determined.\\
It follows from the definition of $A$ and from (\ref{w_h_reln_1}) \& 9\ref{w_h_reln_2}),
\begin{eqnarray}
\int_{-1}^0 W_{20}(\theta)\; d\eta (\theta)&=& 2i\omega_0\tau_0 W_{20}(0)-H_{20}(0),\label{w20_int}\\
\int_{-1}^0 W_{11}(\theta)\; d\eta (\theta)&=& -H_{11}(0),\label{w11_int},
\end{eqnarray}
where $\eta(\theta)=\eta(0,\theta)$.
Also, for $\theta=0$, we have from (\ref{g_fun_main}), (\ref{h_series}) \& (\ref{h_series_g}),
\begin{equation}\label{h20_g}
H_{20}(0)=-g_{20}q(0)-\bar{g}_{02}\bar{q}(0)+2\tau_0(M_1,M_2)^T,
\end{equation}
where $$M_1=m_{13}+m_{14}\gamma,\; M_2=m_{23}\gamma^2+(m_{24}+m_{26}\gamma)\e^{-2i\omega_0\tau_0}+(m_{25}+m_{27}\gamma)\gamma\e^{-i\omega_0\tau_0},$$
and similarly,
\begin{equation}\label{h11_g}
H_{11}(0)=-g_{11}q(0)-\bar{g}_{11}\bar{q}(0)+2\tau_0(N_1,N_2)^T,
\end{equation}
where $$N_1=m_{13}+m_{14}\bar{\gamma},\; M_2=m_{23}\bar{\gamma}^2+(m_{24}+m_{26}\bar{\gamma})\e^{2i\omega_0\tau_0}+(m_{25}+m_{27}\bar{\gamma})\bar{\gamma}\e^{i\omega_0\tau_0}.$$
Noticing that
\begin{eqnarray}
\left(i\omega_0\tau_0 I-\int_{-1}^0\e^{i\omega_0\tau_0}\; d\eta(\theta)\right)q(0)&=& 0,\label{identity1}\\
\left(-i\omega_0\tau_0 I-\int_{-1}^0\e^{-i\omega_0\tau_0}\; d\eta(\theta)\right)\bar{q}(0)&=& 0.\label{identity2}
\end{eqnarray}
Replacing $W_{20}$ in (\ref{w20_int}) using (\ref{w20_soln}),
\begin{equation}\label{e1_int}
\left(2i\omega_0\tau_0 I-\int_{-1}^0\e^{2i\omega_0\tau_0}\; d\eta(\theta)\right)E_1 = 2\tau_0(M_1,M_2)^T,
\end{equation}
which produces the relation
\begin{equation}\label{e1_eqn}
\left(2i\omega_0\tau_0-\tau_0B-\tau_0G\e^{-2i\omega_0\tau_0}\right)E_1=2\tau_0(M_1,M_2)^T,
\end{equation}
that is
\begin{eqnarray}\label{e1_soln}
\begin{pmatrix}
2i\omega_0-m_{11} & -m_{12}\\
-m_{21}\e^{-2\omega_0\tau_0} & 2i\omega_0-m_{22}+m_{22}\e^{-2i\omega_0\tau_0}
\end{pmatrix}
\begin{pmatrix}
E_1^{(1)}\\
E_1^{(2)}
\end{pmatrix} =2 \begin{pmatrix}
M_1\\
M_2
\end{pmatrix}\nonumber\\
\implies E_1^{(1)}=\frac{2}{\Delta}\bigl(2i\omega_0-m_{2}+m_{22}\e^{-2i\omega_0\tau_0}\bigr)M_1+\frac{2}{\Delta}m_{12}M_2\nonumber\\
\& \quad E_1^{(2)}=\frac{2}{\Delta}\bigl(m_{21}\e^{-2i\omega_0\tau_0}M_1+(2i\omega_0-m
_{11})M_2\bigr),
\end{eqnarray}
where $\Delta=(2i\omega_0-m_{11})\left(2i\omega_0-m_{22}+m_{22}\e^{-2i\omega_0\tau_0}\right)-m_{12}m_{21}\e^{-2i\omega_0\tau_0}$.\\
Proceeding in the same way, replacing $W_{20}$ in (\ref{w20_int}) using (\ref{w20_soln}),
\begin{equation}\label{e2_int}
\left(\int_{-1}^0d\eta(\theta)\right)E_2 = -\tau_0(N_1,N_2)^T,
\end{equation}
which produces the relation
\begin{equation}\label{e2_eqn}
(B+G)E_2=-(N_1,N_2)^T,
\end{equation}
that is
\begin{eqnarray}\label{e2_soln}
\begin{pmatrix}
m_{11} & m_{12}\\
m_{21} & 0
\end{pmatrix}
\begin{pmatrix}
E_2^{(1)}\\
E_2^{(2)}
\end{pmatrix} = -\begin{pmatrix}
N_1\\
N_2
\end{pmatrix} \nonumber\\
\implies E_2^{(1)}=-\frac{N_2}{m_{21}} \; \& \; E_2^{(2)}=\frac{m_{11}}{m_{12}m_{21}}N_2-\frac{1}{m_{12}}N_1.
\end{eqnarray}
Now, using (\ref{w20_soln}) \& (\ref{w11_soln}), we can compute $g_{21}$ ad then derive the following values:
\begin{eqnarray}\label{dir_st_hopf_exp}
c_1(0) &=& \frac{i}{2\omega_0\tau_0}\left(g_{20}g_{11}-2|g_{11}|^2-\frac{|g_{02}|^2}{3}\right)+\frac{g_{21}}{2},\nonumber\\
\mu_2 &=& 7 -\frac{\real\{c_1(0)\}}{\real\{\lambda '(\tau_0)\}},\nonumber\\
\beta_2 &=& 2\real\{c_1(0)\},\nonumber\\
T_2 &=& -\frac{\imag\{c_1(0)\}+\mu_2\imag\{\lambda '(\tau_0)\}}{\omega_0\tau_0}
\end{eqnarray}
which give a description of the Hopf bifurcation of system (\ref{sys}) for $\tau=\tau_0$ on the center manifold \& information about the periodic solutions formed for $\tau >\tau_0$.
\begin{theorem}
For the expressions given in (\ref{dir_st_hopf_exp}), the following results hold:
\begin{enumerate}
\item The sign of $\mu_2$ determines the direction of the Hopf bifurcation. If $\mu_2 > 0$, then the Hopf bifurcation is supercritical, and the bifurcating periodic solutions exist for $\tau > \tau_0$. If $\mu_2 < 0$, then the Hopf bifurcation is subcritical, and the bifurcating periodic solutions exist for $\tau < \tau_0$;
\item The parameter $\beta_2$ determines the stability of the bifurcating periodic solutions. The bifurcating periodic solutions are stable if $\beta_2 < 0$ and unstable if $\beta_2 > 0$;
\item Also, $T_2$ determines the period of the bifurcating periodic solutions. The period of the bifurcating periodic solutions increases if $T_2 > 0$ and decreases if $T_2 < 0$.
\end{enumerate}
\end{theorem}

\section{Numerical simulation} \label{numericalsec}
In this section, we present some numerical results of system (\ref{sys}) for a set of parameters and different delay values. The following numerical simulations have been carried out in support of the theoretical results obtained and for a clear illustration of the characteristics and behaviour of the model. All the theoretical findings are verified by the numerical simulations.\\
The parameters chosen for simulation are:
\begin{center}
\begin{tabular}{|c|c|c|c|c|c|c|c|c|}
\hline
Parameter & R & K & M & p & C & D & E & A\\
\hline
Value & 1.8 & 2 & 1.4 & 2 & 10 & 0.2 & 0.5 & 1\\
\hline
\end{tabular}
\end{center}
After non-dimensionalisation of the variables and the parameters, the new parameter set is:
\begin{center}
\begin{tabular}{|c|c|c|c|c|c|c|}
\hline
Parameter & m & p & c & d & e & a\\
\hline
Value & 0.3889 & 2 & 2.5 & 0.2222 & 0.2778 & 0.5\\
\hline
\end{tabular}
\end{center}

Considering the above parameter set, we can compute the interior equilibrium point $(x^*,y^*)$ as $(0.75,1.9688)$. From the theoretical results of the study of the Hopf bifurcation from section \ref{hopfsec} we obtained that Hopf bifurcation occurs at the interior equilibrium point when the $\tau$ value reaches the critical value $\tau_0=1.7805$. Also, from section \ref{dir_st_hopfsec}, using the expressions in (\ref{dir_st_hopf_exp}), we can obtain the type of Hop bifurcation along with the stability nature and the period of the periodic orbits which are formed when $\tau$ value crosses the critical value $\tau_0$.
\begin{figure}[H]
 \subfloat[\label{subfig-7}]{%
  \includegraphics[width=0.5\textwidth, height=6cm]{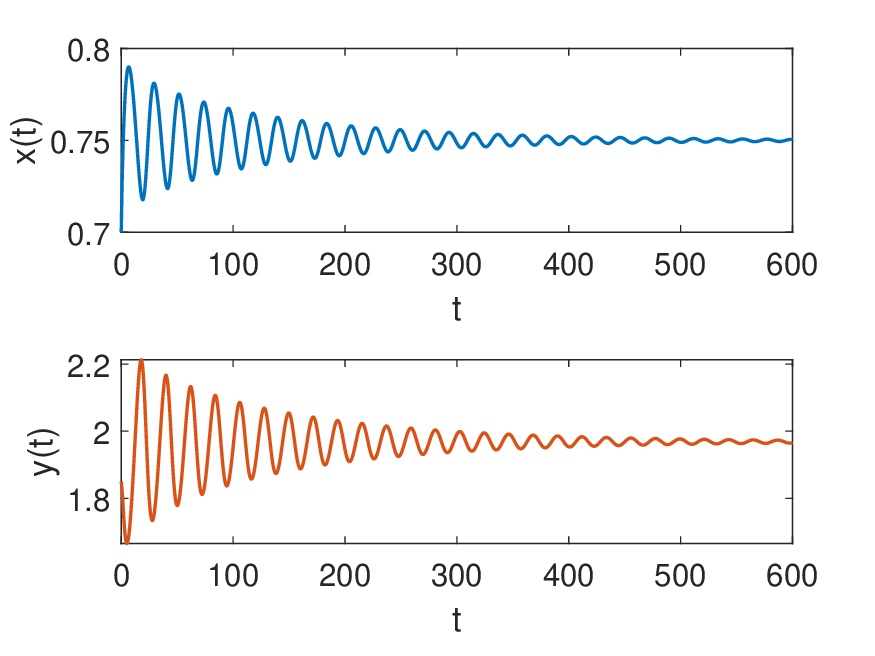}}
 \subfloat[\label{subfig-7}]{%
  \includegraphics[width=0.5\textwidth, height=6cm]{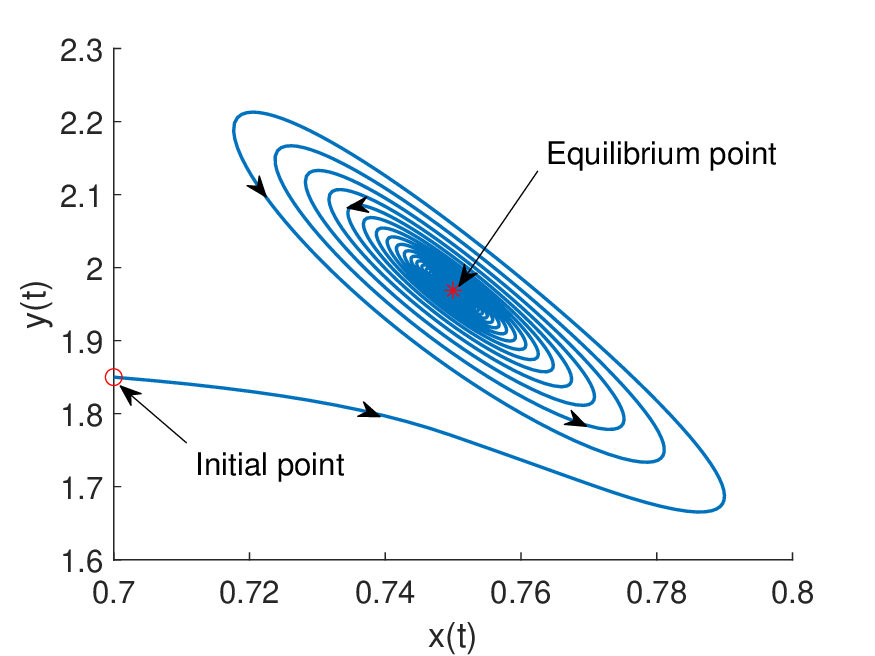}}
 \caption{(a) The time series plots \& (b) the phase portrait of system (\ref{sys}) for delay value $\tau<\tau_0$.}
 \label{bb}
\end{figure}
The above figure shows the time series(on the left) and the phase portrait(on the right) of the populations for the above mentioned parameter set with delay value $\tau=1.75<\tau_0$. We observe that the interior equilibrium point is locally asymptotically stable, so the populations converge to the equilibrium state as time increases.
\begin{figure}[H]
 \subfloat[\label{subfig-7}]{%
  \includegraphics[width=0.5\textwidth, height=6cm]{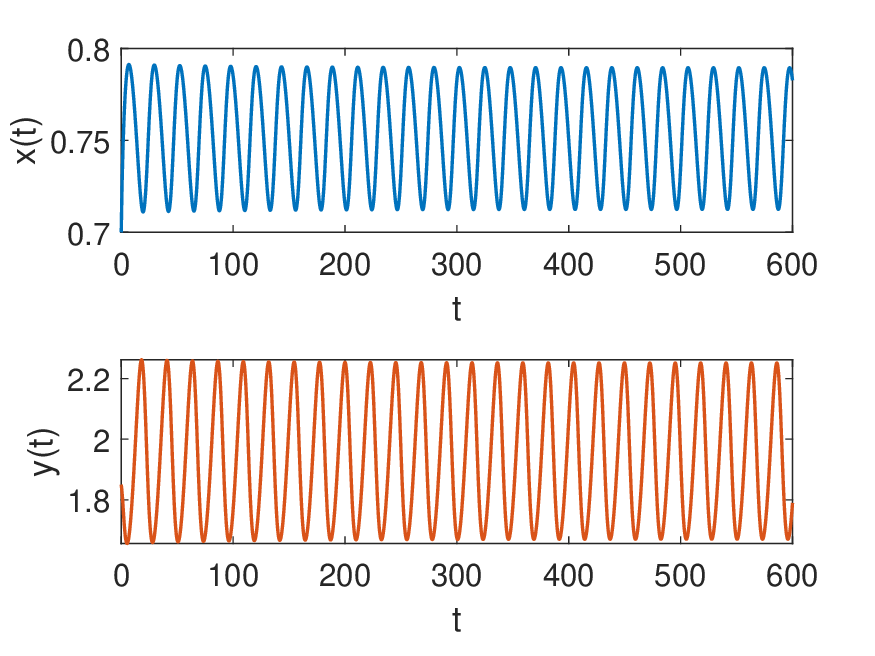}}
 \subfloat[\label{subfig-7}]{%
  \includegraphics[width=0.5\textwidth, height=6cm]{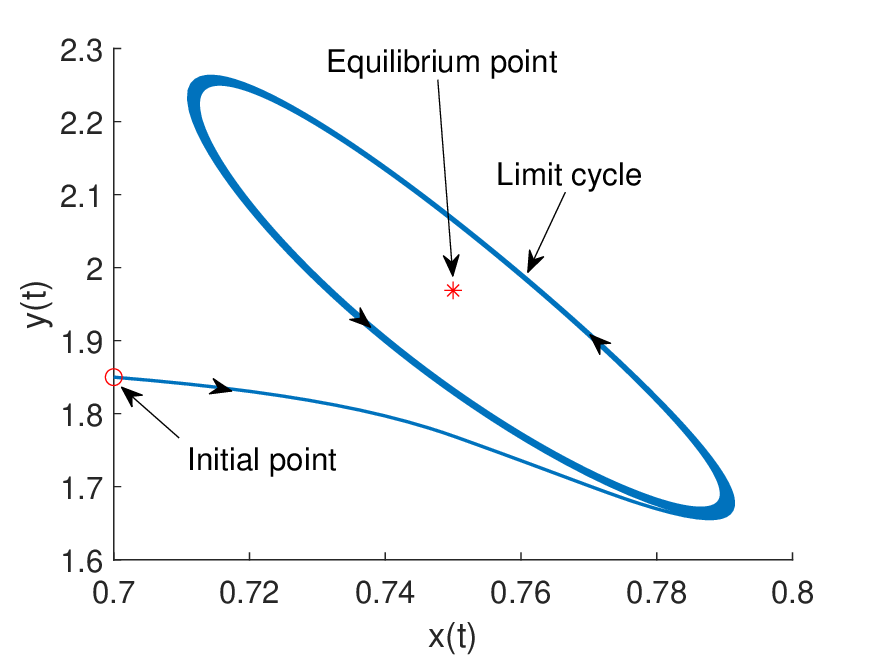}}
 \caption{(a) The time series plots
  \& (b) the phase portrait of system (\ref{sys}) for delay value $\tau>\tau_0$.}
 \label{bb}
\end{figure}
The above figure shows the time series(in the left) and the phase portrait(in the right) of the populations for the above mentioned parameter set with delay value $\tau=1.8>\tau_0$. We observe that the interior equilibrium point is locally asymptotically unstable, and the populations converge to a stable state of oscillation as time increases.

\begin{figure}[H]
 \subfloat[stability region in $a-c$ plane\label{par_space1}]{%
  \includegraphics[width=0.5\textwidth, height=6cm]{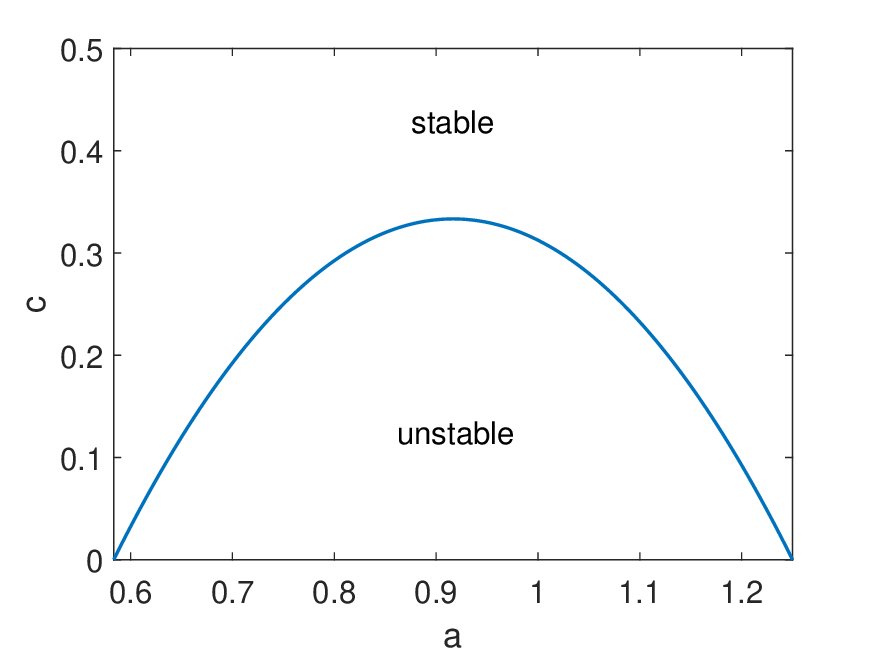}}
 \subfloat[stability region in $c-d$ plane\label{par_space2}]{%
  \includegraphics[width=0.5\textwidth, height=6cm]{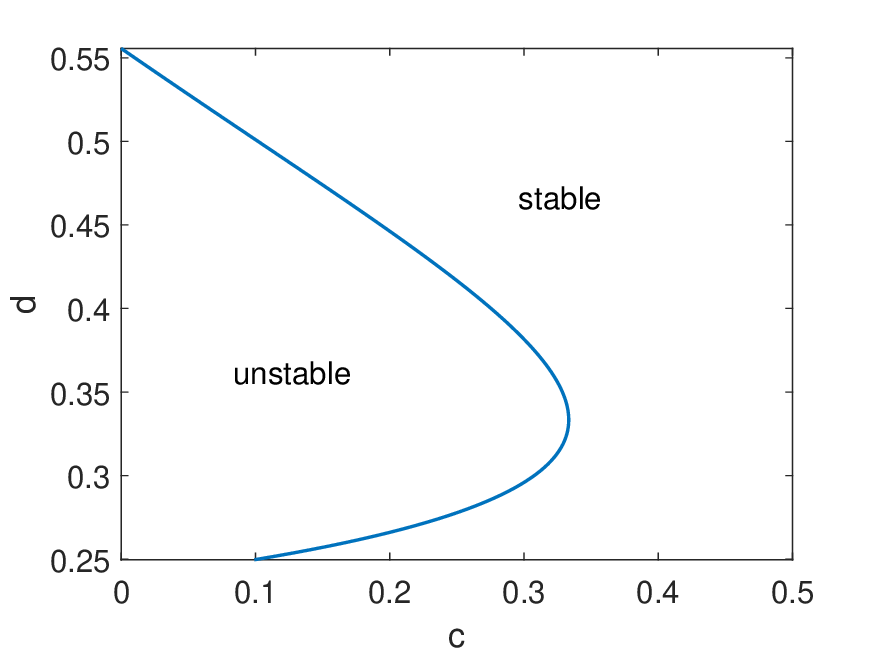}}\\
 \subfloat[stability region in $e-c$ plane\label{par_space3}]{%
  \includegraphics[width=0.5\textwidth, height=6cm]{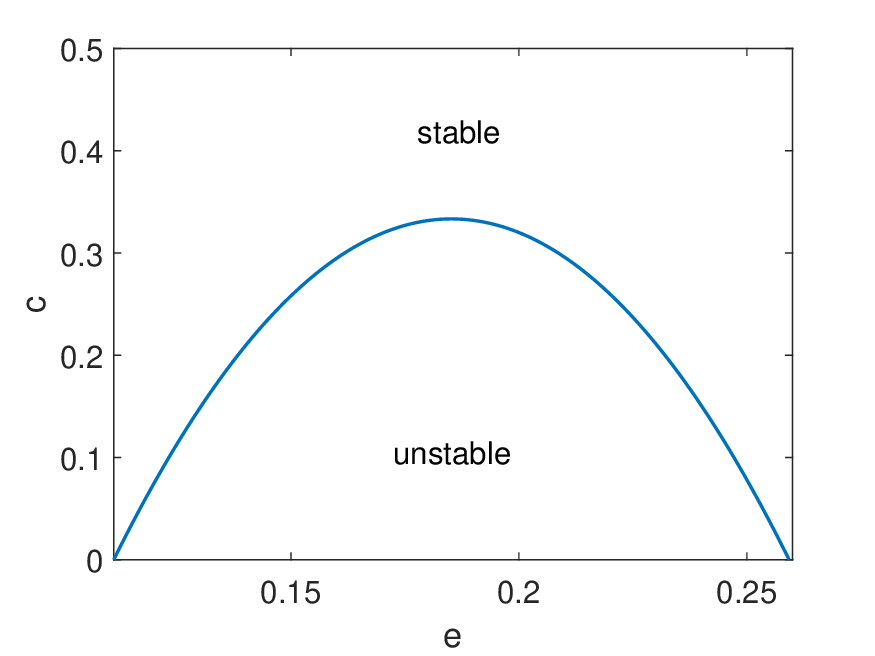}}
 \caption{The graphs depict the local asymptotic stability of the interior equilibrium point as two parameters vary simultaneously.}
 \label{par_space}
\end{figure}
The above graph shows the stability region of the equilibrium point $(x^*,y^*)$ in the parameter space, which is observed during the local stability analysis of system (\ref{sys}). In figure \ref{par_space1} \& \ref{par_space3}, when the respective parameters take the values that lie below the curve, then the interior equilibrium state is locally asymptotically unstable, and when the values lie in the region above the curve, the equilibrium point is locally asymptotically stable. Similarly, in figure \ref{par_space2}, if the parameters take values that lie in the region left to the curve, then the equilibrium point is unstable, and when the parameters take values from the region that lie right to the curve, then the equilibrium point is locally asymptotically stable.
  
\begin{figure}[H]
 \subfloat[PRCC of the parameters for the prey species\label{prcc1}]{%
  \includegraphics[width=0.5\textwidth, height=6cm]{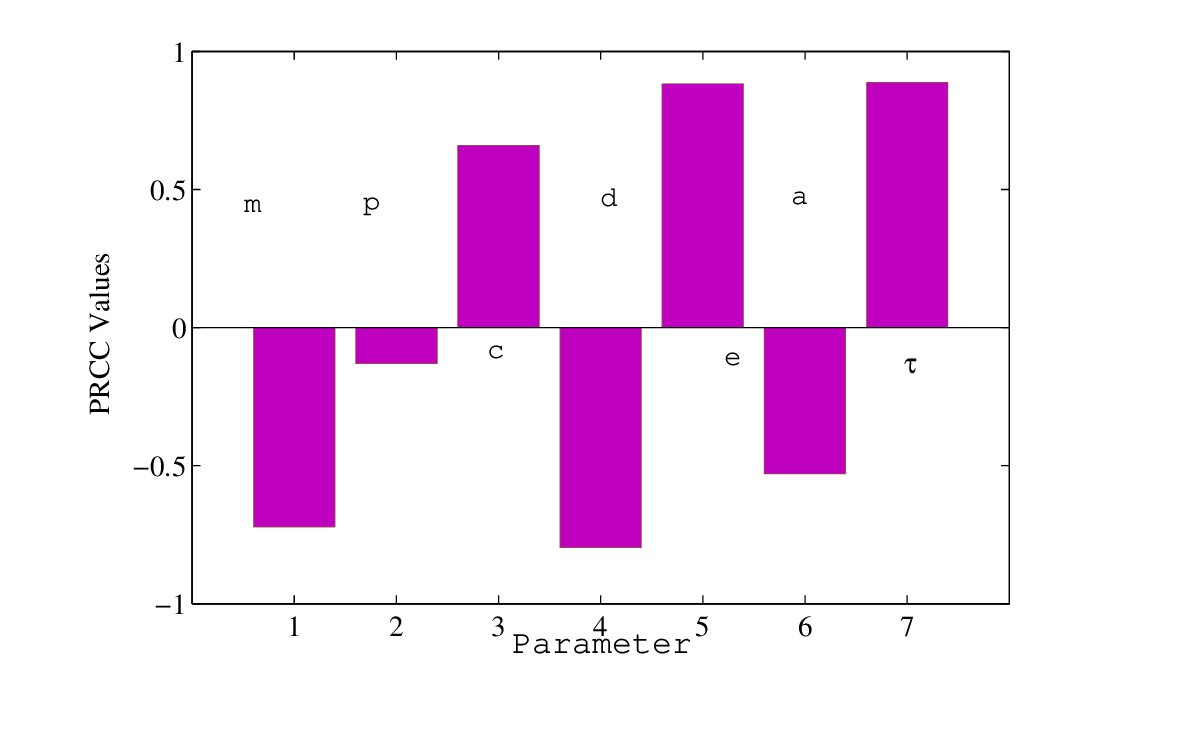}}
 \subfloat[PRCC of the parameters for the predator species\label{prcc2}]{%
  \includegraphics[width=0.5\textwidth, height=6cm]{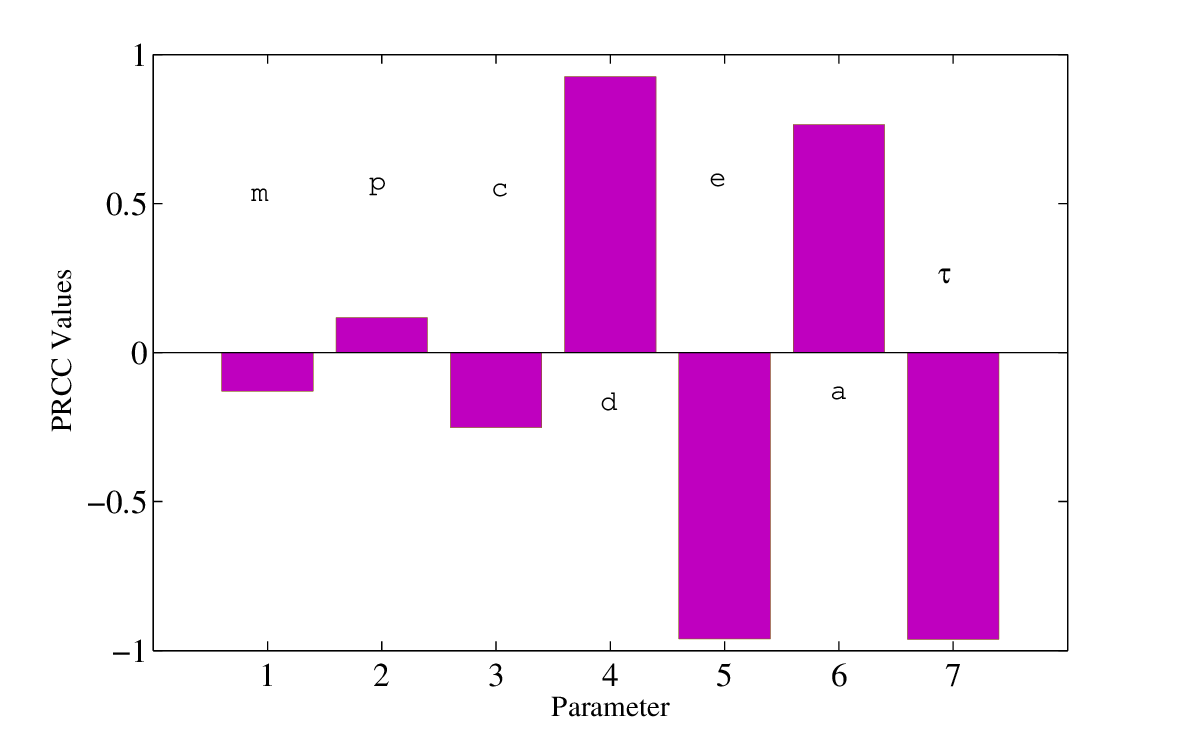}}\\
 \subfloat[PRCC has been plotted over time\label{prcc_ts}]{%
  \includegraphics[width=0.8\textwidth, height=7cm]{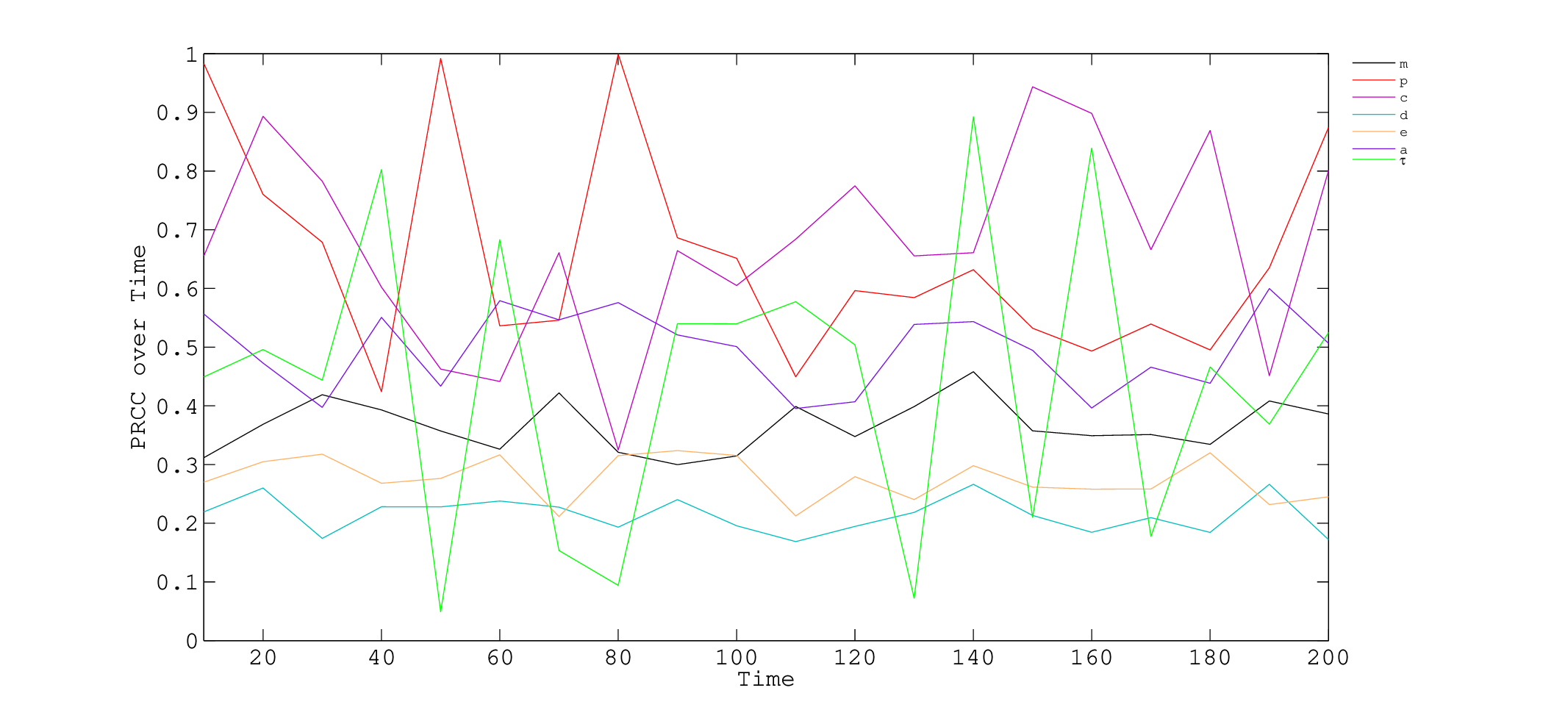}}
 \caption{In this figure, PRCC of different parameters has been plotted for system (\ref{sys}).}
 \label{prcc}
\end{figure}
While analysing a population model, it is very beneficial to investigate the sensitivity of the model to manipulation of the system parameters. This helps us understand which parameter can be considered fundamental for this process\cite{bianca2012} and also, we can explore the level of uncertainty in the parameters. As a result, we can have a good estimate of the parameter value, i.e., parameters those affect more to the model output should have assigned accurate values, whereas those which affect less to the model output can be assigned a rough estimated value\cite{mckay2000}.

Here, we have used a methodology called partial rank correlation coefficient(PRCC). This is the most reliable statistical approach for an effective sensitivity analysis that uses Latin-Hypercube-Sampling(LHS). The PRCC values lie within a range $-1$ to $1$. During the analysis, uniform dissemination is assigned to each parameter of(\ref{sys}) and sampling is done autonomously. The sign of the PRCC values provides the nature of the correlation between the model parameters and the outcome, whereas the magnitude shows the strength of the linear association. We have obtained the PRCC values for the delayed model given in (\ref{sys}) and the results are shown in figure \ref{prcc}.

\section{Population blow-up and its control} \label{controlsec}
In this section, we studied the finite time population explosion and a method to control such situations using biological control. We consider the corresponding non-delayed model system.
\subsection{Sufficient condition for population blow-up}
From system (\ref{sys0}), we observe that, $X'(t)\leq RX(1-X/K)$, so $X(t)\leq K,\;\forall t\geq T_1>0$. So, the population $X$ will never explode.\\
We can say that, $Y(t)$ explodes if
\begin{equation*}
D-\frac{E}{X+A}>0, \quad \forall t\geq 0,
\end{equation*}
which is true if the following condition holds
\begin{equation}\label{blowup_cond}
D-\frac{E}{A}>0.
\end{equation}
Thus we got the sufficient condition for finite time explosion of $Y$ population.\\
Also, from (\ref{blowup_cond}) \& the second equation of (\ref{sys0}),
\begin{equation*}
\frac{dY}{dT} >\left(D-\frac{E}{A}\right)Y^2>0,
\end{equation*}
which gives
\begin{equation*}
Y(t)>\frac{Y_0}{1-Y_0\left(D-\frac{E}{A}\right)t},
\end{equation*}
so, if $T_b$ is the blow-up time then, $$T_b\leq \frac{1}{\left(D-\frac{E}{A}\right)Y_0},$$ where $Y_0$ is the initial population of $Y(T)$.\\
\begin{remark}
Since $D-E/A>0$, so $E/D-A<0$ and the system will not possess any positive interior equilibrium point.
\end{remark}
\subsection{Controlling blow-up using Z-type dynamic method}
In this section, we have used $Z$-type dynamic method\cite{zhang2016} to control the blow-up situation in the system. This is a very efficient control method which forces the dynamics to converge to a desired state by making the error functions vanish, where the error function is the difference of the current state of the variable is supposed to be controlled and the desired state of the variable. A great benefit of using this control is that the rate of convergence can be predetermined so that the convergence can be achieved within a desired time limit.

Let $U_{prey}(t)$ and $U_{pred}(t)$ are two control functions. Incorporating them in the system, we have
\begin{eqnarray}
\frac{dX}{dt}&=& RX\left(1-\frac{X}{K}\right)-\frac{MXY}{X^p+C}-U_{prey}(t)X,\nonumber\\
\frac{dY}{dt}&=& \left(D-\frac{E}{X+A}\right)Y^2-U_{pred}(t)Y.
\end{eqnarray}
These are are to be determined so that the solutions of the system do not blow-up.\\
Let us define the error functions $e_i(t)$, for $i=1,2$, by
$$e_1(t)=X(t)-X_d(t),\quad \text{and}\quad e_2(t)=Y(t)-Y_d(t),$$
such that $$\lim_{t\rightarrow\infty} X(t)=X_d(t) \quad\text{and}\quad \lim_{t\rightarrow\infty} Y(t)=Y_d(t).$$
Let us define the dynamics of the error functions in such a way that the functions exponentially converge to zero at a fixed rate, say $\lambda$. The the dynamics of the error functions are given by
\begin{eqnarray}\label{sys_ctrl}
\dot{e}_1(t)&=& -\lambda e_1(t),\nonumber\\
\dot{e}_2(t)&=& -\lambda e_2(t),
\end{eqnarray}
and thus the controlled system is given by
\begin{eqnarray}
\dot{X}(t)&=& -\lambda(X(t)-X_d(t))+\dot{X}_d(t),\nonumber\\
\dot{Y}(t)&=& -\lambda(Y(t)-Y_d(t))+\dot{Y}_d(t).
\end{eqnarray}
From the above system, it is easy to find out that,
\begin{eqnarray}
U_{prey}(t)&=& R\left(1-\frac{X}{K}\right)-\frac{MY}{X^p+C}+\frac{1}{X}\left[\lambda(X-X_d)-\dot{X}_d(t)\right],\nonumber\\
U_{pred}(t)&=& \left(D-\frac{E}{X+A}\right)Y+\frac{1}{Y}\left[\lambda(Y-Y_d)-\dot{Y}_d(t)\right].
\end{eqnarray}
{\bf Stability:} To ensure the global stability of the only equilibrium point $(0,0)$ of system (\ref{sys_ctrl}), let us construct a positive semi-definite function $V(t)$ such that
\begin{equation}
V(e_1(t),e_2(t))=\frac{1}{2}(e_1^2+e_2^2),
\end{equation}
where $V(0,0)=0$ and $V(e_1,e_2)>0$ if $(e_1,e_2)\neq (0,0)$.\\
Calculating the time derivative of function $V$ along the solutions of (\ref{sys_ctrl}), we have,
\begin{equation}
\dot{V}=-\lambda(e_1^2+e_2^2)<0, \;\;\forall (e_1,e_2)\neq (0,0).
\end{equation}
Hence, $V$ is a Lyapunov function about the equilibrium point $(0,0)$ of (\ref{sys_ctrl}) and $(0,0)$ is globally asymptotically stable. Therefore, $X(t)\rightarrow X_d(t)$ and $Y(t)\rightarrow Y_d(t)$ as $t\rightarrow\infty $.

\begin{figure}[H]
 \subfloat[Time series for prey density\label{subfig-7}]{%
  \includegraphics[width=0.45\textwidth, height=6cm]{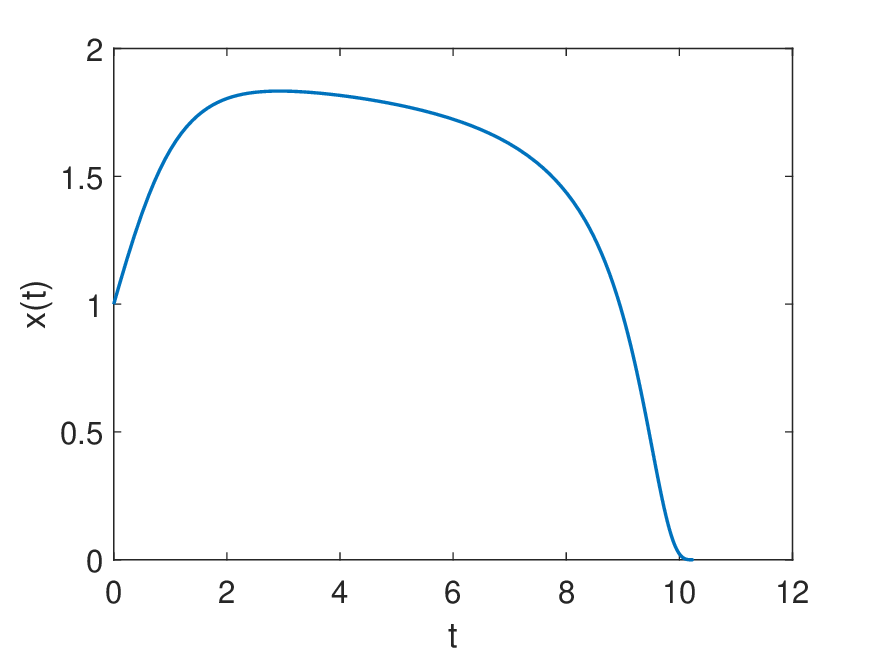}}
 \subfloat[Time series for predator density\label{subfig-7}]{%
  \includegraphics[width=0.45\textwidth, height=6cm]{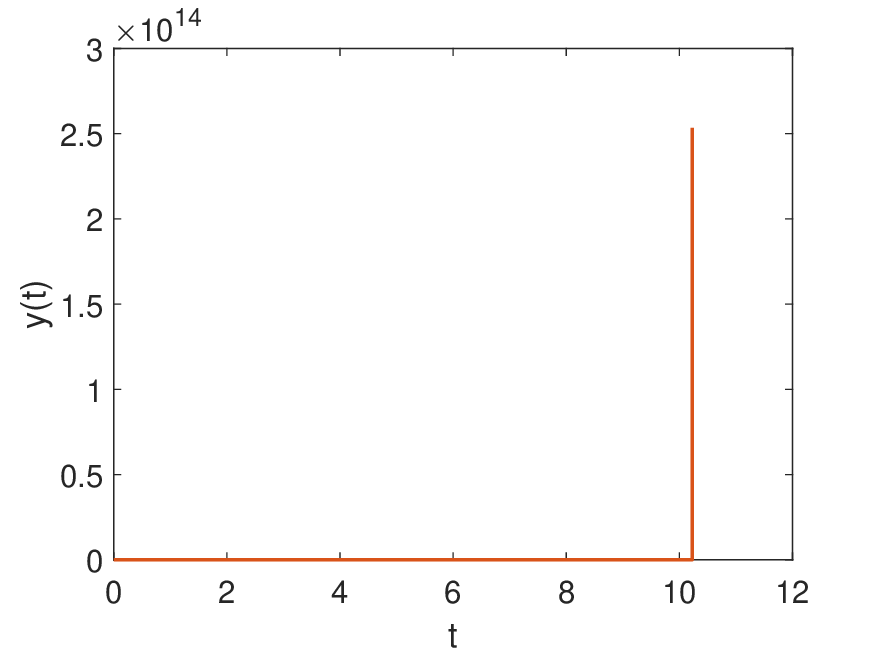}}
 \caption{$Y$ blows up for paramaters $R=1.8$, $K=2$, $M=1.4$, $p=2$, $C=10$, $D=0.2$, $E=0.4$ \& $A=2.5$ satisfying (\ref{blowup_cond}).}
 \label{aa}
\end{figure}
\begin{figure}[H]
 \subfloat[Time series plots for $x(t)$ \& $y(t)$\label{subfig-7}]{%
  \includegraphics[width=0.45\textwidth, height=6cm]{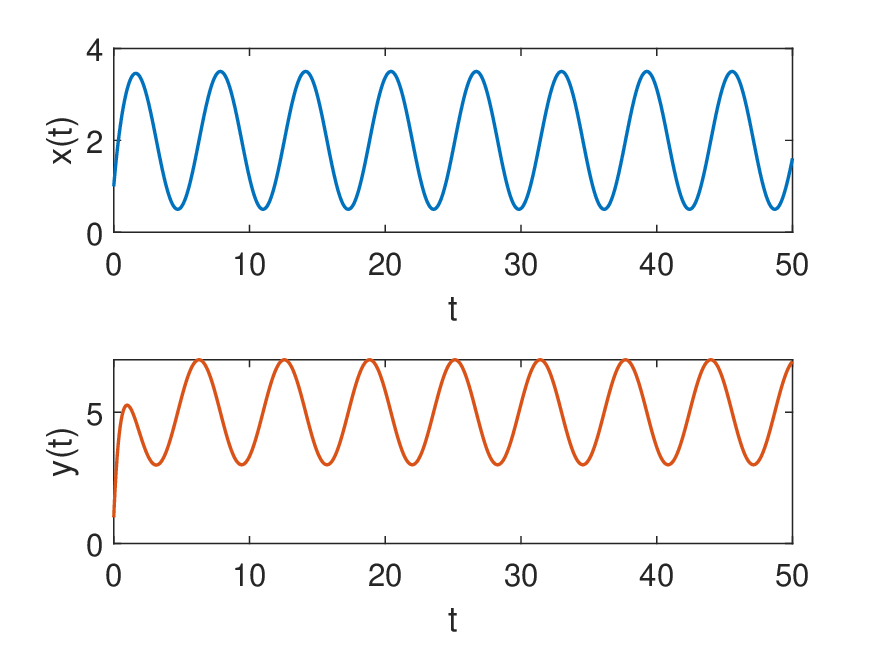}}
 \subfloat[Phase portrait\label{subfig-7}]{%
  \includegraphics[width=0.45\textwidth, height=6cm]{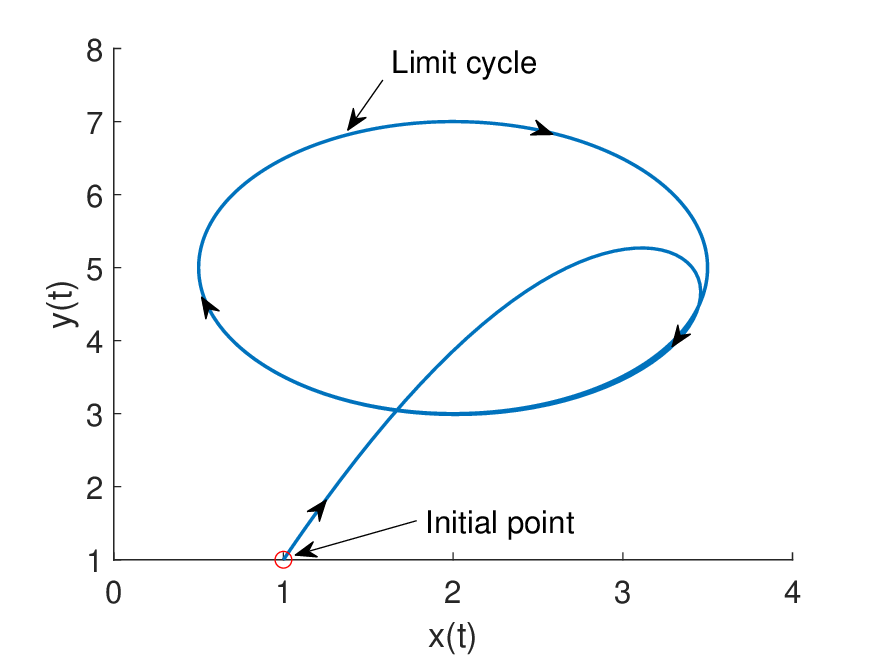}}
 \caption{Figures showing populations achieving desired dynamics when $Z$-type dynamic method is imposed to control the blow-up of $Y$ species. System parameters are given by $R=1.8$, $K=2$, $M=1.4$, $p=2$, $C=10$, $D=0.2$, $E=0.4$ \& $A=2.5$.}
 \label{bb}
\end{figure}
  
\section{Conclusion} \label{conclusionsec}
In this paper, we have considered a system with a prey species showing herd behaviour and a generalist predator species with gestation delay. Due to the food switching behaviour of the generalist predators, it has been a topic of interest for researchers in both field studies and theoretical analysis for many decades\cite{hanski1991,snyder1999,quevedo2009,wimp2019}. Such behaviour has a stabilizing effect on the prey populations as the less numbered species are spared from predation\cite{jaworski2013}. As pointed out by Symondson\cite{symondson2002}, generalist predators can be effective biocontrol agents. Field studies show that they can lessen the impact of crop-destroying harmful pests to a significant degree. To include such characteristics, a modified Leslie-Gower type model is considered to show the predator dynamics where the predator can grow in the absence of the prey species and also, there are additional food sources available for the predators.

The local asymptotic stability of the interior equilibrium point is shown by constructing a suitable Lyapunov function and the stability conditions are obtained. Also, it is observed that Hopf bifurcation occurs at the interior equilibrium point when the delay parameter $\tau$ crosses $\tau_0$. So, when $\tau$ crosses $\tau_0$, a stable limit cycle appears out of the stable equilibrium point and the equilibrium point becomes unstable. The condition for Hopf bifurcation is also obtained. Thus when the delay parameter crosses a threshold value, the population densities, which were going to stabilize at constant densities, in the long run, will now show oscillations over time due to the emergence of the limit cycles.

From the PRCC analysis, we observed that the prey population is most sensitive and positively correlated to the gestation delay value($\tau$) and the death rate of the predators($e$); also, the predator population is most sensitive to these parameters, but in this case, there is a negative correlation.

It is seen that the predator population blows up for a certain parametric restriction. A sufficient condition is obtained for the predator population to blow up. When such a situation arises, the dynamic can be controlled and be made convergent to the desired population dynamic using the $Z-$type dynamic method. Two control functions are obtained using this method, which are functions of the prey density, predator density and the desired population dynamics that can be incorporated in the population dynamics, separately to the prey and predator equations, to obtain the desired dynamics.
\section*{Funding:} Not applicable.
\section*{Data Availability Statement:} Not applicable.
\section*{Declarations:}
\subsection*{Conflict of interest:} The authors declare that they have no conflict of interest.
\subsection*{Code availability:} Not applicable.






\end{document}